\documentclass{amsproc}%
\usepackage{graphicx}
\usepackage{amscd}
\usepackage{amsmath}
\usepackage{amsfonts}
\usepackage{amssymb}%
\theoremstyle{plain}

\numberwithin{equation}{section}

\begin{document}
\title[$\mathbb{Z}_{2}$-graded Cayley-Hamilton trace identities in $\mathrm{M}%
_{n}(E)$]{$\mathbb{Z}_{2}$-graded Cayley-Hamilton trace identities in $\mathrm{M}%
_{n}(E)$}
\author{Szilvia Homolya}
\address{Institute of Mathematics, University of Miskolc, \\
3515 Miskolc-Egyetemv\'{a}ros, Hungary}
\email{szilvia.homolya@gmail.com, mathszil@uni-miskolc.hu}
\author{Jen\H{o} Szigeti}
\address{Institute of Mathematics, University of Miskolc, \\
3515 Miskolc-Egyetemv\'{a}ros, Hungary}
\email{matjeno@uni-miskolc.hu}
\thanks{The second named author was partially supported by the National Research,
Development and Innovation Office of Hungary (NKFIH) K119934 and K135103.}
\thanks{}
\subjclass{15A15,15A24,15A75,15B33,16R40}
\keywords{the full matrix algebra over the infinite dimensional Grassmann algebra, the
Cayley-Hamilton identity}

\begin{abstract}
We show, how the combination of the Cayley-Hamilton theorem and a certain
companion matrix construction can be used to derive $\mathbb{Z}_{2}$-graded
trace identities in $\mathrm{M}_{n}(E)$.

\end{abstract}
\maketitle

\noindent1. INTRODUCTION

\bigskip

\noindent Throughout the paper an algebra $R$\ means a not necessarily
commutative unitary algebra over a commutative ring $C$ (or over a field $K$).
The notation for the full $n\times n$ matrix algebra over $R$ is
$\mathrm{M}_{n}(R)$.

\noindent In case of $\mathrm{char}(K)=0$, Kemer's pioneering work (see
[K1],[K2]) on the $T$-ideals of associative $K$-algebras (leading to the
solution of the Specht problem) revealed the importance of the identities
satisfied by $\mathrm{M}_{n}(E)$ and $\mathrm{M}_{n,d}(E)$, where
$E=E_{0}\oplus E_{1}$ is the naturally $\mathbb{Z}_{2}$-graded Grassmann
(exterior) algebra (over $K$) generated by the infinite sequence of
anticommutative indeterminates $(v_{i})_{i\geq1}$. We note that $E$ is Lie
nilpotent of index $2$.

Let $K\left\langle x_{1},x_{2},\ldots,x_{i},\ldots\right\rangle $ denote the
polynomial $K$-algebra generated by the infinite sequence $x_{1},x_{2}%
,\ldots,x_{i},\ldots$\ of non-commuting indeterminates. The prime $T$-ideals
of this (free associative $K$-) algebra\ are exactly the $T$-ideals of the
identities satisfied by $\mathrm{M}_{n}(K)$ for $n\geq1$ (see [A]). The
$T$-prime (or verbally prime) $T$-ideals are the prime $T$-ideals plus the
$T$-ideals of the identities of $\mathrm{M}_{n}(E)$ for $n\geq1$ and of
$\mathrm{M}_{n,d}(E)$ for $1\leq d\leq n-1$, where $\mathrm{M}_{n,d}(E)$ is
the $K$-subalgebra of the so-called $(n,d)$ supermatrices in $\mathrm{M}%
_{n}(E)$. Another remarkable result is that for a sufficiently large $n\geq1$,
any $T$-ideal contains the $T$-ideal of the identities satisfied by
$\mathrm{M}_{n}(E)$.

The above mentioned three classes of $T$-prime (verbally prime) $PI$-algebras
serve as basic building blocks in Kemer's theory, where $\mathbb{Z}_{2}%
$-graded identities also play an important role. Since the appearance of [K1]
and [K2] considerable efforts have been concentrated on the study of the
various algebraic properties of $\mathrm{M}_{n}(E)$ and $\mathrm{M}_{n,d}(E)$,
see [BR], [DiV], [Do], [KT], [SSz], [Re], [S1], [S2], [V1], [V2].

Accordingly, the importance of matrices over non-commutative rings is an
evidence in the theory of $PI$-rings, nevertheless this fact has been obvious
for a long time in other branches of algebra (structure theory of semisimple
rings, $K$-theory, quantum matrices, etc.). The Cayley-Hamilton theorem and
the corresponding trace identity play a fundamental role in proving classical
results about the polynomial and trace identities of $\mathrm{M}_{n}(K)$. Thus
any Cayley-Hamilton type identity for $\mathrm{M}_{n}(E)$\ seems to be of
general interest.

The main aim of our work is to show, how the combination of the classical
Cayley-Hamilton theorem and a companion matrix construction (from [SSz]) can
be used to derive new $\mathbb{Z}_{2}$-graded trace identities in
$\mathrm{M}_{n}(E)=\mathrm{M}_{n}(E_{0})\oplus\mathrm{M}_{n}(E_{1})$. A
similar approach is used in [HSz] to obtain $\mathbb{Z}_{2}$-graded polynomial
identities in $\mathrm{M}_{n}(E)$.

\bigskip

\noindent2. THE\ APPLICATIONS OF THE CAYLEY-HAMILTON IDENTITY

\bigskip

The Grassmann algebra%
\[
E=K\left\langle v_{1},v_{2},...,v_{i},...\mid v_{i}v_{j}+v_{j}v_{i}=0\text{
for all }1\leq i\leq j\right\rangle =K\left\langle V\right\rangle
\]
over a field $K$\ of characteristic zero (it has been already mentioned in the
introduction) generated by (the countably) infinite set $V=\{v_{1}%
,v_{2},\ldots,v_{t},\ldots\}$ of anticommuting indeterminates can naturally be
extended as%
\[
F=K\left\langle \{w\}\cup V\right\rangle =K\left\langle w,v_{1},v_{2}%
,\ldots,v_{t},\ldots\right\rangle
\]
by using a bigger set $\{w\}\cup V$ of anticommuting generators, where
$w\notin V$. Now we have $v_{i}v_{j}+v_{j}v_{i}=0$, $v_{i}w+wv_{i}=0$ for all
$1\leq i\leq j$ and $w^{2}=0$. Since the cardinalities of $V$ and $\{w\}\cup
V$ are both equal to $\aleph_{0}$, the $K$-algebras $E$ and $F$ are isomorphic.

If $H\in\mathrm{M}_{n}(C)$ is an $n\times n$ matrix over a unitary commutative
ring $C$ with $\frac{1}{m}\in C$ for all integers $m\geq1$, then its
characteristic polynomial $p_{H}(x)=\det(xI-H)\in C[x]$\ can be written as%
\[
p_{H}(x)=\lambda_{0}+\lambda_{1}x+\cdots+\lambda_{n-1}x^{n-1}+\lambda_{n}%
x^{n},
\]
where $\lambda_{n}=1$, $\mathrm{tr}(H)$\ denotes the trace of $H$\ and for
$0\leq k\leq n-1$ we have the following Faddeev-LeVerrier descending recursion
(closely related to\ Newton's identities concerning the elementary symmetric
polynomials)%
\[
\lambda_{k}=-\frac{1}{n-k}\left\{  \lambda_{k+1}\mathrm{tr}(H)+\cdots
+\lambda_{n-1}\mathrm{tr}(H^{n-k-1})+\lambda_{n}\mathrm{tr}(H^{n-k})\right\}
=
\]%
\[
-\frac{1}{n-k}\left\{  \overset{n-k}{\underset{i=1}{\sum}}\lambda
_{k+i}\mathrm{tr}(H^{i})\right\}  .
\]
\noindent\textbf{2.1. Theorem.}\textit{ If }$A\in\mathrm{M}_{n}(E_{0}%
)$\textit{ and }$B\in\mathrm{M}_{n}(E_{1})$\textit{, then a Cayley-Hamilton
trace identity of the form}%
\[
\beta_{0}I_{n}+\underset{k=1}{\overset{n}{\sum}}\left\{  \beta_{k}A^{k}%
+\alpha_{k}(A^{k-1}B+A^{k-2}BA+\cdots+ABA^{k-2}+BA^{k-1})\right\}  =0
\]
\textit{holds, where }$I_{n}\in\mathrm{M}_{n}(E_{0})$\textit{\ is the identity
matrix,}%
\[
p_{A}(x)=\alpha_{0}+\alpha_{1}x+\cdots+\alpha_{n-1}x^{n-1}+\alpha_{n}%
x^{n}=\det(xI-A)\in E_{0}[x]
\]
\textit{is the characteristic polynomial of }$A$\textit{, }$\alpha_{n}=1$,
$\beta_{n}=0$\textit{\ and for }$0\leq k\leq n-1$\textit{ we have}%
\[
\beta_{k}=\left\{  -\frac{1}{n-k}\overset{n-k}{\underset{i=1}{\sum}}%
\beta_{k+i}\mathrm{tr}(A^{i})\right\}  +\left\{  -\frac{1}{n-k}\underset
{r+s\leq n-k-1,0\leq r,s}{\sum}\alpha_{k+r+s+1}\mathrm{tr}(A^{r}%
BA^{s})\right\}  .
\]

\bigskip

\noindent\textbf{Proof.} The so called companion matrix $A+wB$ is in
$\mathrm{M}_{n}(F_{0})$, where $F_{0}$ is the even (and commutative) part of
the extended Grassmann algebra $F=F_{0}\oplus F_{1}$. Any element $\lambda\in
F_{0}$ can be written as $\lambda=\alpha+w\beta$, where $\alpha\in E_{0}$ and
$\beta\in E_{1}$\ are uniquely determined by $\lambda$. Thus we have
$\lambda_{k}=\alpha_{k}+w\beta_{k}$ for the coefficients of the characteristic
polynomial%
\[
p_{A+wB}(x)=\lambda_{0}+\lambda_{1}x+\cdots+\lambda_{n-1}x^{n-1}+\lambda
_{n}x^{n}=
\]%
\[
(\alpha_{0}+w\beta_{0})+(\alpha_{1}+w\beta_{1})x+\cdots+(\alpha_{n-1}%
+w\beta_{n-1})x^{n-1}+(\alpha_{n}+w\beta_{n})x^{n},
\]
where $\alpha_{n}=1$, $\beta_{n}=0$\ and $\alpha_{k}\in E_{0}$, $\beta_{k}\in
E_{1}$ for all $0\leq k\leq n-1$. Using $Aw=wA$, $Bw=-wB$ and $w^{2}=0$, for
the exponent $1\leq i$\ we obtain that%
\[
(A+wB)^{i}=A^{i}+w(A^{i-1}B+A^{i-2}BA+\cdots+ABA^{i-2}+BA^{i-1})
\]
and%
\[
\mathrm{tr}\{\!(A+wB)^{i}\!\}\!=\!\mathrm{tr}(\!A^{i}\!)+w\mathrm{tr}%
(\!A^{i-1}B\!)+w\mathrm{tr}(\!A^{i-2}BA\!)+\cdots+w\mathrm{tr}(\!ABA^{i-2}%
\!)+w\mathrm{tr}(\!BA^{i-1}\!).
\]
In view of the Faddeev-LeVerrier recursion, we deduce that%
\[
\alpha_{k}+w\beta_{k}=-\frac{1}{n-k}\overset{n-k}{\underset{i=1}{\sum}}%
(\alpha_{k+i}+w\beta_{k+i})\{\mathrm{tr}(A^{i})+w\mathrm{tr}(A^{i-1}%
B)+w\mathrm{tr}(A^{i-2}BA)+\cdots
\]%
\[
\cdots+w\mathrm{tr}(ABA^{i-2})+w\mathrm{tr}(BA^{i-1})\}=
\]%
\[
\left\{  -\frac{1}{n-k}\overset{n-k}{\underset{i=1}{\sum}}\alpha
_{k+i}\mathrm{tr}(A^{i})\right\}  +w\left\{  -\frac{1}{n-k}\overset
{n-k}{\underset{i=1}{\sum}}\beta_{k+i}\mathrm{tr}(A^{i})\right\}  +
\]%
\[
+w\left\{  \!-\frac{1}{n-k}\overset{n-k}{\underset{i=1}{\sum}}\alpha
_{k+i}\{\mathrm{tr}(A^{i-1}B)\!+\!\mathrm{tr}(A^{i-2}BA)\!+\cdots
+\!\mathrm{tr}(ABA^{i-2})\!+\!\mathrm{tr}(BA^{i-1})\}\!\right\}  ,
\]
whence%
\[
\alpha_{k}=-\frac{1}{n-k}\overset{n-k}{\underset{i=1}{\sum}}\alpha
_{k+i}\mathrm{tr}(A^{i})
\]
and%
\[
\beta_{k}=\left\{  -\frac{1}{n-k}\overset{n-k}{\underset{i=1}{\sum}}%
\beta_{k+i}\mathrm{tr}(A^{i})\right\}  +\left\{  -\frac{1}{n-k}\underset
{r+s\leq n-k-1,0\leq r,s}{\sum}\alpha_{k+r+s+1}\mathrm{tr}(A^{r}%
BA^{s})\right\}
\]
can be derived. Thus%
\[
p_{A}(x)=\alpha_{0}+\alpha_{1}x+\cdots+\alpha_{n-1}x^{n-1}+\alpha_{n}%
x^{n}=\det(xI-A)\in E_{0}[x]
\]
is the characteristic polynomial of $A\in\mathrm{M}_{n}(E_{0})$.

The application of the Cayley-Hamilton theorem to $A+wB\in\mathrm{M}_{n}%
(F_{0})$ gives that%
\[
0=p_{A+wB}(A+wB)=(\alpha_{0}+w\beta_{0})I_{n}+(\alpha_{1}+w\beta
_{1})(A+wB)+\cdots
\]%
\[
\cdots+(\alpha_{n-1}+w\beta_{n-1})(A+wB)^{n-1}+(\alpha_{n}+w\beta
_{n})(A+wB)^{n}=
\]%
\[
(\alpha_{0}+w\beta_{0})I_{n}+\underset{k=1}{\overset{n}{\sum}}(\alpha
_{k}+w\beta_{k})\{A^{k}+w(A^{k-1}B+A^{k-2}BA+\cdots+ABA^{k-2}+BA^{k-1})\}=
\]%
\[
\left\{  \alpha_{0}I_{n}+\underset{k=1}{\overset{n}{\sum}}\alpha_{k}%
A^{k}\right\}  +w\left\{  \beta_{0}I_{n}+\underset{k=1}{\overset{n}{\sum}%
}\beta_{k}A^{k}\right\}  +
\]%
\[
+w\left\{  \underset{k=1}{\overset{n}{\sum}}\alpha_{k}(A^{k-1}B+A^{k-2}%
BA+\cdots+ABA^{k-2}+BA^{k-1})\right\}  .
\]
Using the fact that for $\alpha\in E_{0}$ and $\beta\in E_{1}$\ the equality
$\alpha+w\beta=0$ implies $\alpha=\beta=0$, we get%
\[
\beta_{0}I_{n}+\underset{k=1}{\overset{n}{\sum}}\left\{  \beta_{k}A^{k}%
+\alpha_{k}(A^{k-1}B+A^{k-2}BA+\cdots+ABA^{k-2}+BA^{k-1})\right\}  =0.\text{
}\square
\]

\noindent\textbf{2.2. Corollary.}\textit{ The case }$n=2$\textit{\ in Theorem
2.1 gives that}%
\[
\left\{  \frac{1}{2}\mathrm{tr}(B)\mathrm{tr}(A)+\frac{1}{2}\mathrm{tr}%
(A)\mathrm{tr}(B)-\frac{1}{2}\mathrm{tr}(AB)-\frac{1}{2}\mathrm{tr}%
(BA)\right\}  I_{2}-\mathrm{tr}(B)A-\mathrm{tr}(A)B+
\]%
\[
AB+BA=0
\]
\textit{holds for }$A\in\mathrm{M}_{2}(E_{0})$\textit{ and }$B\in
\mathrm{M}_{2}(E_{1})$\textit{.}

\bigskip

\noindent\textbf{Proof.} We have%
\[
\beta_{0}I_{2}+\beta_{1}A+\alpha_{1}B+\alpha_{2}(AB+BA)=0
\]
with $\alpha_{2}=1$, $\alpha_{1}=-\mathrm{tr}(A)$, $\beta_{1}=-\mathrm{tr}(B)$
and%
\[
\beta_{0}=-\frac{1}{2}\beta_{1}\mathrm{tr}(A)+\left(  -\frac{1}{2}\right)
\left\{  \alpha_{1}\mathrm{tr}(B)+\alpha_{2}\mathrm{tr}(AB)+\alpha
_{2}\mathrm{tr}(BA)\right\}  =
\]%
\[
\frac{1}{2}\mathrm{tr}(B)\mathrm{tr}(A)+\frac{1}{2}\mathrm{tr}(A)\mathrm{tr}%
(B)-\frac{1}{2}\mathrm{tr}(AB)-\frac{1}{2}\mathrm{tr}(BA).\text{ }\square
\]

\bigskip

\noindent\textbf{2.3. Theorem.}\textit{ If }$B\in\mathrm{M}_{n}(E_{1}%
)$\textit{, then}%
\[
\delta_{0}I_{n}+\gamma_{1}B+\delta_{1}B^{2}+2\gamma_{2}B^{3}+\delta_{2}%
B^{4}+\cdots+(n-1)\gamma_{n-1}B^{2n-3}+\delta_{n-1}B^{2n-2}+nB^{2n-1}=0
\]
\textit{is a Cayley-Hamilton trace identity of degree }$2n-1$\textit{ with
invertible leading coefficient and}%
\[
p_{B^{2}}(x)=\gamma_{0}+\gamma_{1}x+\cdots+\gamma_{n-1}x^{n-1}+\gamma_{n}%
x^{n}=\det(xI-B^{2})\in E_{0}[x]
\]
\textit{is the characteristic polynomial of }$B^{2}\in\mathrm{M}_{n}(E_{0}%
)$\textit{, }$\gamma_{n}=1$\textit{, }$\delta_{n}=0$\textit{\ and for}

\noindent$0\leq k\leq n-1$\textit{ we have}%
\[
\delta_{k}\!=\!\left\{  \!-\frac{1}{n-k}\overset{n-k}{\underset{i=1}{\sum}%
}\delta_{k+i}\mathrm{tr}(B^{2i})\!\right\}  \!+\!\left\{  \!-\frac{1}%
{n-k}\underset{r+s\leq n-k-1,0\leq r,s}{\sum}\gamma_{k+r+s+1}\mathrm{tr}%
(B^{2r+2s+1})\!\right\}  .
\]

\bigskip

\noindent\textbf{Proof.} The substitution $A=B^{2}\in\mathrm{M}_{n}(E_{0})$ in
Theorem 2.1 gives that%
\[
0\!=\!\delta_{0}I_{n}\!+\underset{k=1}{\overset{n}{\sum}}\left\{  \delta
_{k}B^{2k}\!+\!\gamma_{k}(B^{2k-2}B\!+\!B^{2k-4}BB^{2}\!+\cdots+\!B^{2}%
BB^{2k-4}\!+\!BB^{2k-2})\right\}  \!=
\]%
\[
\delta_{0}I_{n}+\underset{k=1}{\overset{n}{\sum}}\left\{  k\gamma_{k}%
B^{2k-1}+\delta_{k}B^{2k}\right\}  ,
\]
where
\[
p_{B^{2}}(x)=\gamma_{0}+\gamma_{1}x+\cdots+\gamma_{n-1}x^{n-1}+\gamma_{n}%
x^{n}=\det(xI-B^{2})\in E_{0}[x]
\]
is the characteristic polynomial of $B^{2}$, $\gamma_{n}=1$, $\delta_{n}%
=0$\ and for $0\leq k\leq n-1$ we have%
\[
\delta_{k}\!=\!\left\{  \!-\frac{1}{n-k}\overset{n-k}{\underset{i=1}{\sum}%
}\delta_{k+i}\mathrm{tr}(B^{2i})\!\right\}  \!+\!\left\{  \!-\frac{1}%
{n-k}\underset{r+s\leq n-k-1,0\leq r,s}{\sum}\gamma_{k+r+s+1}\mathrm{tr}%
(B^{2r+2s+1})\!\right\}  \!.\text{ }\square
\]

\bigskip

\noindent\textbf{2.4. Remark.} The direct application of the Cayley-Hamilton
theorem to $B^{2}\in\mathrm{M}_{n}(E_{0})$ gives a monic identity for $B$\ of
degree $2n$ with coefficients in $E_{0}$. The degree of our (also can be
considered as monic) identity in Theorem 2.3 is only $2n-1$, where the even
degree coefficients are in $E_{1}$ and the odd degree coefficients are in
$E_{0}$.

\bigskip

\noindent\textbf{2.5. Corollary.}\textit{ If }$B\in\mathrm{M}_{n}(E_{1}%
)$\textit{, }$2\leq n$\textit{\ and }$\mathrm{tr}(B^{t})=0$\textit{ for all
}$1\leq t\leq2n-2$\textit{, then }$nB^{2n-1}=\mathrm{tr}(B^{2n-1})I_{n}%
$\textit{ is a scalar matrix. The additional condition }$\mathrm{tr}%
(B^{2n-1})=0$\textit{\ implies that }$B^{2n-1}=0$\textit{.}

\bigskip

\noindent\textbf{Proof.} In view of the recursions ($\gamma_{n}=1$,
$\delta_{n}=0$, $0\leq k\leq n-1$)%
\[
\gamma_{k}=-\frac{1}{n-k}\overset{n-k}{\underset{i=1}{\sum}}\gamma
_{k+i}\mathrm{tr}(B^{2i})
\]
and%
\[
\delta_{k}\!=\!\left\{  \!-\frac{1}{n-k}\overset{n-k}{\underset{i=1}{\sum}%
}\delta_{k+i}\mathrm{tr}(B^{2i})\!\right\}  \!+\!\left\{  \!-\frac{1}%
{n-k}\underset{r+s\leq n-k-1,0\leq r,s}{\sum}\gamma_{k+r+s+1}\mathrm{tr}%
(B^{2r+2s+1})\!\right\}  ,
\]
we obtain that $\gamma_{1}=\cdots=\gamma_{n-1}=0$, $\delta_{1}=\cdots
=\delta_{n-1}=0$ and%
\[
\delta_{0}\!=\!-\frac{1}{n}\underset{r+s\leq n-1,0\leq r,s}{\sum}%
\gamma_{r+s+1}\mathrm{tr}(B^{2r+2s+1})\!=\!-\frac{1}{n}\underset{r+s=n-1,0\leq
r,s}{\sum}\gamma_{r+s+1}\mathrm{tr}(B^{2r+2s+1})\!=
\]%
\[
-\frac{1}{n}\underset{r+s=n-1,0\leq r,s}{\sum}\gamma_{n}\mathrm{tr}%
(B^{2n-1})=-\mathrm{tr}(B^{2n-1})
\]
in Theorem 2.3, whence $nB^{2n-1}=\mathrm{tr}(B^{2n-1})I_{n}$ follows.
$\square$

\bigskip

\noindent\textbf{2.6. Remark.} In Rosset's proof of the Amitsur-Levitzki
theorem ([Ro]) and in [KSz] a similar implication was used, namely that
$\mathrm{tr}(B^{2})=\mathrm{tr}(B^{4})=\cdots=\mathrm{tr}(B^{2n})=0$ implies
$B^{2n}=0$.

\bigskip

\noindent\textbf{2.7. Corollary.}\textit{ The case }$n=2$\textit{\ in Theorem
2.3 gives that}%
\[
\left\{  \mathrm{tr}(B)\mathrm{tr}(B^{2})-\mathrm{tr}(B^{3})\right\}
I_{2}-\mathrm{tr}(B^{2})B-\mathrm{tr}(B)B^{2}+2B^{3}=0.
\]
\textit{holds for }$B\in\mathrm{M}_{2}(E_{1})$\textit{. The case }%
$n=3$\textit{\ in Theorem 2.3 gives that}%
\[
\left\{  -\frac{1}{2}\mathrm{tr}^{2}(B^{2})\mathrm{tr}(B)+\mathrm{tr}%
(B^{3})\mathrm{tr}(B^{2})+\frac{1}{2}\mathrm{tr}(B^{4})\mathrm{tr}%
(B)-\mathrm{tr}(B^{5})\right\}  I_{3}+
\]%
\[
\left\{  \!\frac{1}{2}\mathrm{tr}^{2}(B^{2})\!-\!\frac{1}{2}\mathrm{tr}%
(B^{4})\!\right\}  \!B+\left\{  \!\mathrm{tr}(B^{2})\mathrm{tr}%
(B)\!-\!\mathrm{tr}(B^{3})\!\right\}  \!B^{2}-2\mathrm{tr}(B^{2}%
)B^{3}-\mathrm{tr}(B)B^{4}+3B^{5}\!=\!0
\]
\textit{holds for }$B\in\mathrm{M}_{3}(E_{1})$\textit{.}

\bigskip

\noindent\textbf{Proof.} In case $n=2$ we have%
\[
\delta_{0}I_{2}+\gamma_{1}B+\delta_{1}B^{2}+2B^{3}=0,
\]
where%
\[
p_{B^{2}}(x)=\gamma_{0}+\gamma_{1}x+\gamma_{2}x^{2}=\det(xI-B^{2})\in
E_{0}[x]
\]
is the characteristic polynomial of $B^{2}\in\mathrm{M}_{2}(E_{0})$ with
$\gamma_{2}=1$ and $\gamma_{1}=-\mathrm{tr}(B^{2})$. The recursion
($\delta_{2}=0$, $0\leq k\leq1$)%
\[
\delta_{k}=\left\{  -\frac{1}{2-k}\overset{2-k}{\underset{i=1}{\sum}}%
\delta_{k+i}\mathrm{tr}(B^{2i})\right\}  +\left\{  -\frac{1}{2-k}%
\underset{r+s\leq2-k-1,0\leq r,s}{\sum}\gamma_{k+r+s+1}\mathrm{tr}%
(B^{2r+2s+1})\right\}
\]
gives that $\delta_{1}=-\mathrm{tr}(B)$ and%
\[
\delta_{0}=-\frac{1}{2}\delta_{1}\mathrm{tr}(B^{2})+\left(  -\frac{1}%
{2}\right)  \left\{  \gamma_{1}\mathrm{tr}(B)+\gamma_{2}\mathrm{tr}%
(B^{3})+\gamma_{2}\mathrm{tr}(B^{3})\right\}  =\mathrm{tr}(B)\mathrm{tr}%
(B^{2})-\mathrm{tr}(B^{3}).
\]
In case $n=3$ we have%
\[
\delta_{0}I_{3}+\gamma_{1}B+\delta_{1}B^{2}+2\gamma_{2}B^{3}+\delta_{2}%
B^{4}+3B^{5}=0,
\]
where%
\[
p_{B^{2}}(x)=\gamma_{0}+\gamma_{1}x+\gamma_{2}x^{2}+\gamma_{3}x^{3}%
=\det(xI-B^{2})\in E_{0}[x]
\]
is the characteristic polynomial of $B^{2}\in\mathrm{M}_{3}(E_{0})$. In view
of the recursion ($\gamma_{3}=1$, $0\leq k\leq2$)%
\[
\gamma_{k}=-\frac{1}{3-k}\overset{3-k}{\underset{i=1}{\sum}}\gamma
_{k+i}\mathrm{tr}(B^{2i}),
\]
we obtain that $\gamma_{2}=-\mathrm{tr}(B^{2})$ and%
\[
\gamma_{1}=-\frac{1}{2}\overset{2}{\underset{i=1}{\sum}}\gamma_{1+i}%
\mathrm{tr}(B^{2i})=-\frac{1}{2}\left\{  -\mathrm{tr}^{2}(B^{2})+\mathrm{tr}%
(B^{4})\right\}  .
\]
The other recursion ($\delta_{3}=0$, $0\leq k\leq2$)%
\[
\delta_{k}=\left\{  -\frac{1}{3-k}\overset{3-k}{\underset{i=1}{\sum}}%
\delta_{k+i}\mathrm{tr}(B^{2i})\right\}  +\left\{  -\frac{1}{3-k}%
\underset{r+s\leq3-k-1,0\leq r,s}{\sum}\gamma_{k+r+s+1}\mathrm{tr}%
(B^{2r+2s+1})\right\}
\]
gives that $\delta_{2}=-\mathrm{tr}(B)$,%
\[
\delta_{1}=\left\{  -\frac{1}{2}\delta_{2}\mathrm{tr}(B^{2})\right\}
+\left\{  -\frac{1}{2}\underset{r+s\leq1,0\leq r,s}{\sum}\gamma_{r+s+2}%
\mathrm{tr}(B^{2r+2s+1})\right\}  =
\]%
\[
\frac{1}{2}\mathrm{tr}(B)\mathrm{tr}(B^{2})+\left\{  -\frac{1}{2}\gamma
_{2}\mathrm{tr}(B)-\frac{1}{2}\gamma_{3}\mathrm{tr}(B^{3})-\frac{1}{2}%
\gamma_{3}\mathrm{tr}(B^{3})\right\}  =\mathrm{tr}(B^{2})\mathrm{tr}%
(B)-\mathrm{tr}(B^{3})
\]%
\[
\delta_{0}=\left\{  -\frac{1}{3}\overset{3}{\underset{i=1}{\sum}}\delta
_{i}\mathrm{tr}(B^{2i})\right\}  +\left\{  -\frac{1}{3}\underset
{r+s\leq2,0\leq r,s}{\sum}\gamma_{r+s+1}\mathrm{tr}(B^{2r+2s+1})\right\}  =
\]%
\[
\left\{  -\frac{1}{3}\delta_{1}\mathrm{tr}(B^{2})-\frac{1}{3}\delta
_{2}\mathrm{tr}(B^{4})\right\}  +\left(  -\frac{1}{3}\right)  \left\{
\gamma_{1}\mathrm{tr}(B)+2\gamma_{2}\mathrm{tr}(B^{3})+3\gamma_{3}%
\mathrm{tr}(B^{5})\right\}  =
\]%
\[
\left\{  -\frac{1}{3}\left\{  \mathrm{tr}(B^{2})\mathrm{tr}(B)-\mathrm{tr}%
(B^{3})\right\}  \mathrm{tr}(B^{2})+\frac{1}{3}\mathrm{tr}(B)\mathrm{tr}%
(B^{4})\right\}  +
\]%
\[
\left(  -\frac{1}{3}\right)  \left\{  -\frac{1}{2}\left\{  -\mathrm{tr}%
^{2}(B^{2})+\mathrm{tr}(B^{4})\right\}  \mathrm{tr}(B)-2\mathrm{tr}%
(B^{2})\mathrm{tr}(B^{3})+3\mathrm{tr}(B^{5})\right\}  =
\]%
\[
-\frac{1}{2}\mathrm{tr}^{2}(B^{2})\mathrm{tr}(B)+\mathrm{tr}(B^{3}%
)\mathrm{tr}(B^{2})+\frac{1}{2}\mathrm{tr}(B^{4})\mathrm{tr}(B)-\mathrm{tr}%
(B^{5})
\]
Thus we have%
\[
\delta_{0}I_{3}+\gamma_{1}B+\delta_{1}B^{2}+2\gamma_{2}B^{3}+\delta_{2}%
B^{4}+3B^{5}=
\]%
\[
\left\{  -\frac{1}{2}\mathrm{tr}^{2}(B^{2})\mathrm{tr}(B)+\mathrm{tr}%
(B^{3})\mathrm{tr}(B^{2})+\frac{1}{2}\mathrm{tr}(B^{4})\mathrm{tr}%
(B)-\mathrm{tr}(B^{5})\right\}  I_{3}+
\]%
\[
\left\{  -\frac{1}{2}\left\{  -\mathrm{tr}^{2}(B^{2})+\mathrm{tr}%
(B^{4})\right\}  \right\}  B+
\]%
\[
\left\{  \mathrm{tr}(B^{2})\mathrm{tr}(B)-\mathrm{tr}(B^{3})\right\}
B^{2}-2\mathrm{tr}(B^{2})B^{3}-\mathrm{tr}(B)B^{4}+3B^{5}.\text{ }\square
\]

\bigskip

\noindent REFERENCES

\bigskip

\noindent\lbrack A] S.A. Amitsur, \textit{The T-ideals of the free ring}, J.
London Math. Soc. 30 (1955), 470-475.

\noindent\lbrack BR] A. Berele, A. Regev, \textit{Asymptotic codimensions of
}$\mathrm{M}_{k}(E)$, Advances in Mathematics 363 (2020), 106979.

\noindent\lbrack DiV] O. M. Di Vincenzo, \textit{On the graded identities of
}$\mathrm{M}_{1,1}(E)$, Israel J. Math. 80(3) (1992), 323--335.

\noindent\lbrack Do] M. Domokos, \textit{Cayley-Hamilton theorem for }%
$2\times2$\textit{\ matrices over the Grassmann algebra}, J. Pure Appl.
Algebra 133 (1998), 69-81.

\noindent\lbrack HSz] Sz. Homolya, J. Szigeti, $\mathbb{Z}_{2}$\textit{-graded
identities of }$\mathrm{M}_{n}(E)$\textit{, two descendants of the
Amitsur-Levitzki theorem}, submitted

\noindent\lbrack KT] I. Kantor, I. Trishin, \textit{On a concept of
determinant in the supercase,} Communications in Algebra, 22(10) (1994), 3679--3739.

\noindent\lbrack K1] A.R. Kemer, \textit{Varieties of }$\mathbb{Z}_{2}%
$\textit{-graded algebras}, Math. USSR Izv. 25 (1985), 359-374.

\noindent\lbrack K2] A.R. Kemer, \textit{Ideals of Identities of Associative
Algebras}, Translations of Math. Monographs, Vol. 87 (1991), AMS Providence,
Rhode Island.

\noindent\lbrack KSz] P. K\"{o}rtesi, J. Szigeti, \textit{The adjacency matrix
of a directed graph over the Grassmann algebra}, Algebra and its Applications,
Contemporary Mathematics Vol. 259 (2000), 319-322. (AMS Providence, Rhode Island)

\noindent\lbrack Re] A. Regev, \textit{Tensor products of matrix algebras over
the Grassmann algebra}, J. Algebra 133(2) (1990), 512--526.

\noindent\lbrack Ro] S. Rosset,\textit{ A new proof of the Amitsur-Levitzki
identity}, Israel J. Math. 23 (1976), 187--188.

\noindent\lbrack SSz] S. Sehgal, J. Szigeti, \textit{Matrices over centrally
}$\mathbb{Z}_{2}$\textit{-graded rings}, Beitrage zur Algebra und Geometrie,
Vol. 43 (2002), No.2, 399-406.

\noindent\lbrack S1] J. Szigeti, \textit{New determinants and the
Cayley-Hamilton theorem for matrices over Lie nilpotent rings}, Proc. Amer.
Math. Soc. 125 (1997), 2245-2254.

\noindent\lbrack S2] J. Szigeti, \textit{On the characteristic polynomial of
supermatrices}, Israel J. Math. 107 (1998), 229-235.

\noindent\lbrack V1] U. Vishne, \textit{Polynomial Identities of }$2\times
2$\textit{ matrices over the Grassmannian}, Communications in Algebra, Vol.
30(1) (2002), 443-454.

\noindent\lbrack V2] U. Vishne, \textit{Polynomial Identities of }%
$\mathrm{M}_{2,1}(G)$, Communications in Algebra, Vol. 39(6) (2011), 2044-2050.

\end{document}